\documentclass[12pt]{article}
\usepackage{mathrsfs}
\usepackage{amsfonts}
\usepackage{bbm}
\usepackage{amssymb}
\usepackage{amsmath,amssymb}
\usepackage{latexcad}
\newcommand{\G}{\Gamma}
\newcommand{\nin}{\noindent}
\newcommand{\vs}{\vspace*}

\newcommand{\eset}{\emptyset}

\newcommand{\seq}{\subseteq}

\newcommand{\ol}{\overline}

\def \nin {\noindent}

\def \Lemma #1 {\vs{3mm}\nin {\bf Lemma #1} \it}
\def \Prop #1 {\vs{3mm}\nin {\bf Proposition #1} \it}
\def \Th #1 {\vs{3mm}\nin {\bf Theorem #1} \it}
\def \Cor #1 {\vs{3mm}\nin {\bf Corollary #1} \it}
\def \Ex #1 {\vs{3mm}\nin {\bf Example #1} \it}
\def \Proof {\vs{3mm}\nin {\bf Proof. }}
\def \part #1 {\hfil\break\hglue 12pt {\rm (#1)~}}

\def\fs{\footnotesize}

\title{
\bf\LARGE On a class of semigroup graphs\thanks {This research is
supported by the National Natural Science Foundation of China (Grant
No.10671122). } }
\author{{Li Chen\thanks{chenli\_830202@sjtu.edu.cn}\, and Tongsuo
Wu\thanks{
tswu@sjtu.edu.cn}}\\
 {\fs Department of Mathematics, Shanghai
Jiaotong University, Shanghai 200240, P.R.China}\\
}
\date{}
\begin{document}
\baselineskip=16pt \maketitle

\vs{3mm}\nin{\small\bf Abstract. }

{\fs Let $G=\Gamma(S)$ be a semigroup graph, i.e., a zero-divisor
graph of a semigroup $S$ with zero element $0$. For  any adjacent
vertices $x, y$ in $G$, denote $C(x,y)=\{z\in
V(G)\,|\,N(z)=\{x,y\}\}$. Assume that in $G$  there exist two
adjacent vertices $x,y$,
 a vertex $s\in C(x,y)$ and a vertex $z$
such that $d(s,z)=3$. In this paper, we study algebraic properties
of $S$ with such graphs $G=\G(S)$, giving some sub-semigroups and
ideals of $S$. We  construct some classes of such semigroup graphs
and classify all semigroup graphs with the property in two cases.}

\vs{3mm}\nin {\small Key Words:} {\small Zero-divisor semigroup;
Sub-semigroup; Zero-divisor graph; Graph classification}

\vs{4mm}\nin{\bf 1. Introduction}

\vs{3mm}\nin Throughout, $G$ is a simple and connected graph. For a
vertex $x$ of $G$, the neighborhood of $x$ is denoted as $N(x)$,
which is the set of all vertices adjacent to $x$. Denote also
$\ol{N(x)}=N(x)\cup \{x\}$. The cardinality of $N(x)$ is denoted by
$deg(x)$. The vertex $x$ is called an {\it end vertex} if $deg(x)=
1$, and an {\it isolated vertex} if $deg(x)=0$. Throughout, $S$ is a
commutative semigroup with $0$. Recall that for  a commutative
semigroup (or a commutative ring) $S$ with $0$, the zero-divisor
graph $\Gamma(S)$ is an undirected graph whose vertices are the
zero-divisors of $S\setminus \{0\}$, and with two vertices $a,b$
adjacent in case $ab=0$
(\cite{B},\cite{AL},\cite{DF},\cite{DD},\cite{CW},\cite{WuChen}). If
$G\cong \G(S)$ for some semigroup $S$ with zero element $0$, then
$G$ is called a {\it semigroup graph}.

 Some fundamental
properties and possible algebraic structure of $S$ and graphic
structures of $\Gamma(S)$ were established in \cite{AL, DF, DD}
among others. For example, it was proved that $\Gamma(S)$ is always
connected, and the diameter of $\Gamma(S)$ is less than or equal to
$3$. If $\Gamma(S)$ contains a cycle, then its core, i.e., the union
of the cycles in $\G(S)$, is a union of squares and triangles, and
any vertex not in the core is an end vertex which is connected to
the core by a single edge. In \cite[Theorem 4]{DD}, the authors
provided a descending chain of ideals $I_k$ of $S$, where $I_k$
consists of all elements of $S$ with vertex degree greater than or
equal to $k$ in $\G(S)$. In \cite{WuLu}, the authors continued the
study on the sub-semigroup structure and ideal structure of
semigroups. By \cite[Theorem 1.3]{DMS}, $\G(S)$ contains no cycle if
and only if $\G(S)$ is either a star graph or a two-star graph. By
\cite[Theorem 2.10]{LuWu}, the core $K(G)$ contains no triangle if
and only if $\G(S)$ is a bipartite graph, if and only if $\G(S)$ is
one of the following graphs: star graphs, two-star graphs,  complete
bipartite graphs,   complete bipartite graphs with a thorn. By
\cite[Theorem 2.3]{LiuWu}, the core $K(G)$ contains no rectangle if
and only if $\G(S)$ is one of the following graphs: an isolated
vertex, a star graph, a two-star graph, a triangle with $n$ thorns
($n = 0, 1, 2, 3$), a fan graph, a fan graph with a thorn adjacent
to its center.

Let $S$ be a commutative semigroup with zero-element $0$, and let
$G=\G(S)$. For any adjacent vertices $a, b$ in $V(G)$,  denote
$C(a,b)=\{x\in V(G)\,|\,N(x)=\{a,b\}\}$ and let $T_a$ denote the set
of all end vertices adjacent to $a$. Consider the following
condition assumed on $G=\G(S)$

\vspace{3mm}($\triangle$)  There exist in $G$ two adjacent vertices
$a,b$,
 a vertex $s\in C(a,b)$ and a vertex $z$
such that $d(s,z)=3$.

\vspace{3mm}

\nin In this paper, we study algebraic properties of semigroups $S$
and the graphic structures of $\G(S)$ such that the condition
($\triangle$) holds for $\G(S)$. (We can further assume that
triangles and rectangles coexist in the core $K[\G(S)]$.)  In
particular, it is proved that $S\setminus [C(a, b)\cup T_a\cup T_b]$
is an ideal of $S$ (Theorem 2.4). Under some additional conditions,
it is proved that $S\setminus C(a,b)$ is a sub-semigroup of $S$ and
there exists an element $c\in C(a,b)$ such that $[S\setminus
C(a,b)]\cup \{c\}$ is also a sub-semigroup of $S$. We also use
Theorem 2.4 to construct some classes of semigroup graphs which
satisfies the condition $(\triangle)$, and give a complete
classification of such semigroup graphs in two cases.

We record a known result on finite semigroups to end this part (see,
e.g., \cite[Corollary 5.9 on page 25 ]{Gri}). We also include a
proof for the completeness.

\vs{3mm}\nin {\bf Lemma 1.1.} {\it Any finite nonempty semigroup $S$
contains an idempotent element.}

\Proof Take any element $x$ from $S$ and consider the sequence
$x,x^2,x^3,\cdots.$ Since $S$ is a finite set, there exist $m<n$
such that $x^m=x^n$. Let $r=n-m$, and take $k$ such that $kr\ge m$.
Then$$x^{kr}=x^m\cdot x^{kr-m}=(x^r\cdot x^m) x^{kr-m}=x^r\cdot
x^{kr}=x^r(x^rx^{kr})=\cdots =(x^{kr})^2. \quad \Box$$

\vs{4mm}\nin{\bf 2. Properties of $S$ }

\vs{3mm}\nin{\bf Lemma 2.1.} {\it Let $S$ be a commutative semigroup
with 0, $\Gamma(S)$ its zero-divisor graph. For any vertex $x\in
\Gamma(S)$, if there exists a vertex $y\in \Gamma(S)$ such that
$d(x,y)=3$, then $x^{2}\neq 0$ in $S$. }

\vs{2mm}\nin{\bf Proof.} As $d(x,y)=3$, there exist vertices $a,z\in
\Gamma(S)$ such that $x-a-z-y$,  $xz\neq 0$ and $ ay\neq 0$. If
$x^{2}=0$, then $x^{2}z=0$ and thus $xz\in Ann(x)$. Clearly $xz\in
Ann(y)$. Then $xz\in Ann(x)\cap Ann(y)=\{0\}$, a contradiction.
\quad $\Box$

\vs{2mm} Part of the following result  is contained in
\cite[Proposition 2.8]{WuLu}.

\vs{3mm}\nin{\bf Proposition 2.2.} {\it Let $G=\Gamma(S)$ be a
zero-divisor graph of a semigroup $S$. For a vertex $b\in V(G)$, let
$T_b=\{x\in V(G)|xb=0,\,  x\not=b,\,deg(x)=1\}$.}

(1) {\it If $b^2\neq 0$, then $T_b\cup \{0\}$ is a sub-semigroup of
$S$.}

(2) {\it If $b$ is not an end vertex and $T_b\neq \emptyset$, then
$\{0, b\}$ is an ideal of $S$.}

\vs{2mm}\nin{\bf Proof.} (1) We only need consider as $T_b\neq
\emptyset$. If $G$ contains no cycle, then $G$ is either a two-star
graph or a star graph by \cite[Theorem 1.3]{DMS} . If $G$ is a star
graph, then $T_b=S\setminus \{b,0\}$. For all $x, y\in T_b$, we must
have $xy\neq b$, since otherwise $0=xyb=b^2\neq 0$,  a
contradiction. This shows that $T_b\cup \{0\}$ is a sub-semigroup of
$S$ when $G$ is a star graph. If $G$ is a two-star graph or a graph
with cycles, then $B\neq \emptyset$ where $B=\{x\in
V(G)\,|\,deg(x)\geq 2 , xb=0\}$. For all $x\in T_b$, we have
$$x^2\in Ann(b)\cup \{0\}=\{0\}\cup T_b\cup B$$
\nin If $x^2\in B$, denote $x^2=v$. Then there exists $z\in
S\setminus
 \{b\}$ such that $zv=0$. Since $x^2z=vz=0$, we have $xz\in
Ann(x)=\{x, b\}$. If $xz=x$, then $v=x^2=x^2z=vz=0$, a
contradiction. If $xz=b$, then $0=xzb=b^2\neq 0$,  another
contradiction. So we must have $x^2\in T_b\cup \{0\}$. If $|T_b|\geq
2$, then exists a vertex $y\in T_b$ such that $x\neq y$. If $xy\in
B$, denote $xy=v$. Then there exists $z\in S\setminus \{b\}$ such
that $zv=0$. As $xyz=0$, we have $xz\in Ann(y)=\{y,b\}$ and $yz\in
Ann(x)=\{x,b\}$, and thus $xz=y$ and $yz=x$. Then $x^2=xyz=vz=0$. On
the other hand, $xy=x^2z=0$, a contradiction. So $xy\in T_b\cup
\{0\}$, and hence $T_b\cup \{0\} \leqslant S$.

(2) Since $T_b\neq \emptyset$, there exists $x\in T_b$ such that
$by\in Ann(x)=\{0,b,x\}$ for all $y\in S$. By assumption, $b$ is not
an end vertex and thus there exists $z\in S\backslash \{x\}$ such
that $bz=0$. Then  $by\neq x$ since otherwise, $by=x$ and it implies
$0=bzy=zx\neq 0$, a contradiction. This completes the proof. \quad
$\Box$

\vs{2mm}\nin{\bf Remark 2.3.} {\it In Proposition 2.2(1), the
conclusion can not hold if $b^2=0$. }

\vs{3mm} For a vertex $v$ of a graph $G$, if $v$ is not an end
vertex and there is no end vertex adjacent to $v$, then $v$ is said
to be an {\it internal vertex}.  We know prove the main result of
this section.

\vs{3mm}\nin{\bf Theorem 2.4.} {\it Let $G=\Gamma(S)$ be a semigroup
graph satisfying condition $(\triangle)$. Denote $L=\{z\in S\,|\,d(s,z)=3\}$. Then $\{0,a,b\}$ is an
ideal of $S$,  $S\setminus [C(a, b)\cup T_a\cup T_b]$ is an ideal
of $S$ and $L$ is a sub-semigroup of $S$. Furthermore, }

  (1) {\it If both $a$ and $b$ are internal vertices, then
$S\setminus C(a, b)$ is an ideal of  $S$.}

  (2) {\it If $a$ is an internal vertex, while $b$ is not an internal vertex and $b^2\neq
  0$, then $S\setminus C(a, b)$ is a sub-semigroup of  $S$.}

\vs{2mm}\nin{\bf Proof.} 
Fix some $s\in C(a,b)$ and let $B=\{x\,|\,x\in S,x\notin C(a,b),
\,d(s,x)=2\}$, $L=\{y\,|\,y\in S, \,d(s,y)=3\}$. By assumption
$L\not=\eset$, $C(a,b)\not=\eset$, and $T_a\cup T_b\subset B$.
Notice that there is no end vertex in $B\setminus (T_a\cup T_b)$. By
\cite[Theorem 2.3]{AL} or by \cite[Theorem 1(2)]{DD},
$S=\{0,a,b\}\cup  C(a,b)\cup B\cup L$ and it is a disjoint union of
four nonempty subsets. By Lemma 2.1. we have $c^2\neq 0,\, \forall
c\in C(a,b)$, and hence $Ann(c)=\{0, a, b\}$. Clearly, $a^2\in
Ann(c)$, $b^2\in \{0,a,b\}$ and
$$\{0,a,b\}(B\cup L)\subseteq Ann(c)=\{0,a, b\}.$$
This shows that $\{0,a,b\}$ is an ideal of $S$.

For any $y$ in $L$, there exists a vertex $x\in B$ such that $xy=0$.
Then  $yS\subseteq Ann(x)$ while $C(a,b)\cap Ann(x)=\eset$. Hence
$LS\cap C(a,b)=\eset$. Furthermore, for any $s\in S$, $ys\in
\{0,a,b\}\cup L\cup B$. If $ys\in \{0,a,b\}\cup B$, then it is clear
that $ys\not\in T_a\cup T_b$ whether $sy=x$ or not. Thus $LS\cap
(T_a\cup T_b)=\eset$, and hence $$(C(a,b)\cup T_a\cup T_b)\cap
LS=\eset.$$ For any vertex $x_1$ in $B\setminus (T_a\cup T_b)$,
$x_1\not\in C(a,b)$ and it has degree greater than one. Hence for
any $x_{1}\in B\setminus (T_a\cup T_b)$ and any $x_{2}\in S$, there
exists a vertex $u\in B\cup L $ such that $x_{1}u=0$. Then
$x_{1}x_{2}\in Ann(u)$ and it implies $x_{1}x_{2}\notin C(a,b)$.
Thus $ [(B\setminus (T_a\cup T_b))S]\cap C(a,b)=\eset$. Finally, by
\cite[Theorem 4]{DD}, the core of $G$ together with $0$ forms an
ideal of $S$. Thus these arguments show that $S\setminus [C(a,
b)\cup T_a\cup T_b]$ is an ideal of $S$.

Now take any $u,v\in L$ and consider $uv$. Clearly $uv\not=0$ (see also remark(1) proceeding Theorem 3.3).
Also $uv\not\in C(a,b) $. If $uv\not\in L$, then  we can assume $(uv)a=0$. Since $\{0,a,b\}$ is an ideal of $S$,
we would have either $ua=0$ or $ub=0$. This shows $uv\in L$ and hence $L$ is a sub-semigroup of $S$. Notice $0\not\in L$.

(1) If both $a$ and $b$ are internal vertices, then $S\setminus
C(a,b)=S\setminus [C(a, b)\cup T_a\cup T_b]$. In this case,
$S\setminus C(a, b)$ is clearly an ideal of  $S$.

(2) Now assume that $b$ is not an internal vertex, and $b^2\not=0$.
Again let $T_b$ be the set of end vertices adjacent to $b$. By the
above discussion, we already have $([\{0,a,b\}\cup L\cup (B\setminus
T_b)]S)\cap C(a,b)=\eset$. Since $b^2\neq 0$, we have
$T_b^2\leqslant T_b\cup \{0\}$ by Theorem 2.2(1). These facts show
that $S\setminus C(a, b)$ is a sub-semigroup of $S$, and it
completes the proof. \quad $\Box$

\vs{2mm}\nin{\bf Remarks 2.5.(1)}   In Theorem 2.4,  if there is no
$z\in V(G)$ such that $d(s,z)=3$, then the theorem may not hold. An
example is contained in Example 3.1.

\vs{2mm}\nin{\bf Remarks 2.5.(2)}  Theorem 2.4 can be easily extended for a graph which satisfies the condition $(K_p)$: in the definition $(\triangle)$, replace $a-b$ by an induced complete subgraph $H$ with $p-1$ vertices, for any prime number $p$. The condition $(K_p)$ is rather natural for the ring graphs, see our subsequent work \cite{ChenWu2}.

\vs{3mm}\nin{\bf Theorem 2.6.} {\it Let $G=\Gamma(S)$ be a semigroup
graph satisfying condition $(\triangle)$. If further $a$ is an
internal vertex, while $b$ is not an internal vertex and $b^2\neq
0$, then both $\{0,a\}$ and $\{0,b\}$ are ideals of $S$. }

\vs{2mm}\nin{\bf Proof.}
 By Theorem 2.2(2), we already
have $\{0,b\}\trianglelefteq S$.

Let $B=\{x\,|\,x\notin C(a,b),\,d(c,x)=2\}$, $L=\{y\,|\,d(c,y)=3\}$.
For any $d\in T_b$, we have $ad\in Ann(c)=\{a,b,0\}$. Clearly,
$ad\neq 0$. We conclude $ad=a$. In fact, if $ad\not=a$, then $ad=b$
and hence $ad^2=bd=0$. By Theorem 2.2(1), we get $d^2\in T_b\cup
\{0\}$, and we have $d^2\neq 0$ by Lemma 2.1. Thus $d^2\in T_b$, and
thus $ad^2\neq 0$,
 a contradiction. The contradiction shows $ad=a$.

For any $x\in B\setminus T_b$, we have $ax\in Ann(c)=\{a,b,0\}$. If
$ax=b$, then $b=ax=adx=bd=0$, a contradiction. Thus we must have
$aB\subseteq \{0,a\}$.

In a similar way, we prove $aL\subseteq \{0,a\}$. This completes the
proof. \quad $\Box$

\vs{3mm}In the following we proceed to prove that under some
additional conditions, there exists an element $c$ in $C(a,b)$ such
that $[S\setminus C(a, b)]\cup \{c\}$ is a sub-semigroup of $S$. For
this purpose, we need the following technical lemma.

\vs{3mm}\nin{\bf Lemma 2.7.} {\it Let $G=\Gamma(S)$ be a semigroup
graph satisfying condition $(\triangle)$. Assume further that one of
the following conditions is satisfied:}

(1) {\it Both $a$ and $b$ are internal vertices.}

(2) {\it $a$ is an internal vertex, $b^2\not=0$ and $T_b=\{d\}$.}

\nin {\it Then there exists an element $c\in C(a,b)$ such that
$[S\setminus C(a, b)]\cup \{c\}\leqslant S$, if and only if there
exists an element $c_1\in C(a,b)$ such that $c_1^2\in [S\setminus
C(a,b)]\cup \{c_1\}$.}

\vs{2mm}\nin{\bf Proof.} $\Rightarrow$ Clear. In fact, we further
have $c^2\not=0$ by Lemma 2.1.

$\Leftarrow$ Assume that there exists an element $c_1\in C(a,b)$
such that $c_1^2\in [S\setminus C(a,b)]\cup \{c_1\}$. Then $c_1^2\in
\{a, b, c_1\}\cup X$, where $X$ is the set of the vertices which are
 adjacent to $a,b$ and at the same time belong to $B$.

(1) If  both $a$ and $b$ are internal vertices, then $T_a\cup
T_b=\eset$. In this case, repeat the proof of Theorem 2.4 and obtain
$(c_1B\cup c_1L)\cap C(a,b)=\eset$. Hence $c_1^2=c_1$ implies
$[S\setminus C(a, b)]\cup \{c_1\}\leqslant S$.

(2) Now assume that $a$ is an internal vertex, while $b$ is not an
internal vertex,  $b^2\not=0$ and  $T_b=\{d\}$. By Theorem 2.4(2),
$S\setminus C(a,b)\leqslant S$. Since $S\setminus [C(a,b)\cup T_b]$
is an ideal of $S$, we already have $c_1[S\setminus (C(a,b)\cup
T_b)]\seq S\setminus (C(a,b)$.

If $c_1T_b\seq \{c_1\}\cup [S\setminus C(a,b)]$, then $[S\setminus
C(a, b)]\cup \{c_1\}$ is a sub-semigroup of $S$ since $c_1^2\in
\{c_1\}\cup [S\setminus C(a,b)]$. In the following we assume that
$c_1d\in C(a,b)\setminus \{c_1\}$, and denote $c_1d=c$. Then $d^2=d$
since $c_1d^2=cd\not=0$ and $T_b\cup \{0\}\leqslant S$. Then $cd=c$.
Since $c_1^2\in [S\setminus C(a,b)]\cup \{c_1\}$, we have the
following four possible subcases.

(2.1) $c_1^2=a$. In this case, we have $cd=c$ and
$c_1c=c_1^2d=ad=a$. Then $c^2=c_1cd=ad=a\in Ann(c)=\{a,b\}$ since
$\{0,a,b\}\leqslant S$. Thus $[S\setminus C(a, b)]\cup
\{c\}\leqslant S$.

(2.2) $c_1^2=b$. In this case, $cd=c$ and $cc_1=c_1^2d=bd=0$, a
contradiction. Thus this case can not occur.

(2.3) $c_1^2=c_1$. In this case, $cd=c$ and $cc_1=c_1^2d=c_1d=c$.
Then $c^2=cc_1d=cd=c$, and hence $[S\setminus C(a, b)]\cup
\{c\}\leqslant S$.

(2.4) $c_1^2=x\in X$. In this case, $cd=c$ and $cc_1=c_1^2d=xd\not
\in C(a,b)$ by Theorem 2.4. Then $c^2=c_1cd=xd^2=xd\not \in C(a,b)$,
and hence $[S\setminus C(a, b)]\cup \{c\}\leqslant S$.

This completes the proof. \quad $\Box$

\vs{3mm}Now we are ready to prove

\vs{3mm}\nin{\bf Proposition 2.8.} {\it Let $G=\Gamma(S)$ be a
semigroup graph satisfying condition $(\triangle)$. Assume that
$|C(a,b)|$ is finite. If one of the following conditions is
satisfied, then there exists an element $c\in C(a,b)$ such that
$[S\setminus C(a, b)]\cup \{c\}$ is a sub-semigroup of $S$: }

(1) {\it Both $a$ and $b$ are internal vertices.}

(2) {\it $a$ is an internal vertex, $b^2\not=0$ and $T_b=\{d\}$.}

\vs{2mm}\nin{\bf Proof.} If $C(a,b)$ is a sub-semigroup of $S$, then
 by Lemma 1.1. there is an element $c_1\in C(a,b)$ such that
 $c_1^2=c_1$. By Lemma 2.7, there exists an element $c\in C(a,b)$ such that
$[S\setminus C(a, b)]\cup \{c\}\leqslant S$. In the following we
assume that $C(a,b)$ is not a sub-semigroup of $S$. Then there exist
$c_i, c_j\in C(a,b)$ such that $c_ic_j\not\in C(a,b)$, and this
implies $c_ic_j\in \{a,b\}\cup X$ where $X$ is the set of the
vertices which are adjacent to $a,b$ and at the same time belong to
$B$.

Assume $c_ic_j=a$. If $i=j$, then the result follows from Lemma 2.7.
If $i\neq j$, then $c_ic_j=a$ implies $c_{i}^2c_j=ac_i=0$, and thus
$c_{i}^2\in Ann(c_j)$, i.e. $c_i^2\in \{a, b\}$. Then we use Lemma
2.7 again to obtain the result. When $c_ic_j=b$, a similar
discussion lead to the result.

Finally, assume $c_ic_j=x\in X$. In this case, it is only necessary
to consider the $i\neq j$ case. Since $x\in X$, there is an element
$z\in B\cup L$ such that $xz=0$. Then $c_ic_jz=xz=0$, and hence
$c_jz\in Ann(c_i)=\{0,a,b\}$. Thus $c_j^2z=ac_j=0$ or
$c_j^2z=bc_j=0$, i.e. $c_j^2\in Ann(z)$. This means $c_j^2\in
\{a,b,X\}\cap Ann(z)$. By Lemma 2.7, there exists an element $c\in
C(a,b)$ such that $[S\setminus C(a, b)]\cup \{c\}$ is a
sub-semigroup of $S$. This completes the proof. \quad $\Box$

\vs{0.4cm}\nin{\bf 3. Some Examples and complete classifications of
the graphs in two cases }

\vs{3mm}\nin In this section, we use Theorem 2.4 to study the
correspondence of zero-divisor semigroups and several classes of
graphs satisfying the four necessary conditions of \cite[Theorem
1]{DD} as well as the general assumption of Theorem 2.4.

\vs{2mm}\nin{\bf Example 3.1.}  {\it Consider the graph $G$ in
Fig.3, where both $U$ and $V$ consist of end vertices. We claim that
each graph in Fig.3 is a semigroup graph.}

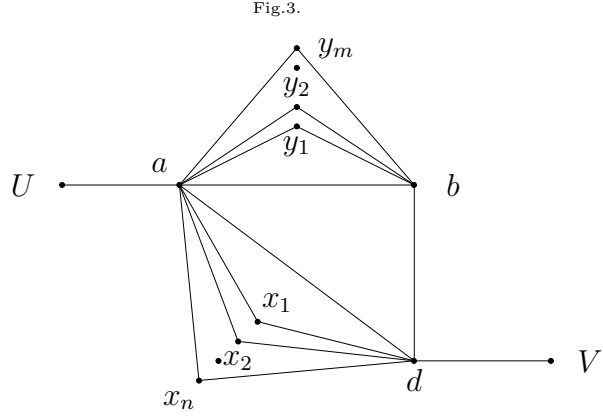
\begin{center}
 \setlength{\unitlength}{0.13cm}
\begin{picture}(72,52)
\thinlines
\drawdot{34.0}{34.0}
\drawdot{22.0}{28.0}
\drawdot{46.0}{28.0}
\drawpath{34.0}{34.0}{22.0}{28.0}
\drawpath{22.0}{28.0}{46.0}{28.0}
\drawpath{46.0}{28.0}{34.0}{34.0}
\drawdot{34.0}{36.0}
\drawdot{34.0}{40.0}
\drawdot{34.0}{42.0}
\drawpath{22.0}{28.0}{34.0}{36.0}
\drawpath{34.0}{36.0}{46.0}{28.0}
\drawpath{46.0}{28.0}{34.0}{42.0}
\drawpath{34.0}{42.0}{22.0}{28.0}
\drawcenteredtext{34.0}{32.0}{$y_1$}
\drawcenteredtext{34.0}{38.0}{$y_2$}
\drawcenteredtext{38.0}{42.0}{$y_m$}
\drawcenteredtext{32.0}{46.0}{\tiny Fig.3.}
\drawdot{46.0}{10.0}
\drawdot{60.0}{10.0}
\drawpath{60.0}{10.0}{46.0}{10.0}
\drawpath{46.0}{10.0}{46.0}{28.0}
\drawpath{22.0}{28.0}{46.0}{10.0}
\drawdot{30.0}{14.0}
\drawdot{28.0}{12.0}
\drawdot{24.0}{8.0}
\drawpath{46.0}{10.0}{30.0}{14.0}
\drawpath{30.0}{14.0}{22.0}{28.0}
\drawpath{22.0}{28.0}{28.0}{12.0}
\drawpath{28.0}{12.0}{46.0}{10.0}
\drawpath{46.0}{10.0}{24.0}{8.0}
\drawpath{24.0}{8.0}{22.0}{28.0}
\drawdot{10.0}{28.0}
\drawcenteredtext{50.0}{28.0}{$b$}
\drawcenteredtext{46.0}{8.0}{$d$}
\drawcenteredtext{32.0}{16.0}{$x_1$}
\drawcenteredtext{28.0}{10.0}{$x_2$}
\drawcenteredtext{22.0}{6.0}{$x_n$}
\drawcenteredtext{64.0}{10.0}{$V$}
\drawpath{10.0}{28.0}{22.0}{28.0}
\drawcenteredtext{6.0}{28.0}{$U$}
\drawcenteredtext{20.0}{30.0}{$a$}
\drawdot{26.0}{10.0}
\end{picture}\par
\end{center}

In fact, first notice that $d(y_i, V)=3,\,C(a,b)=\{y_1,..., y_m\}$,
and $C(a,d)=\{x_1,...,x_n\}$. By Theorem 2.4, if $G$ has a
corresponding semigroup $S=V(G)\cup \{0\}$, then the subset
$S\setminus (\{y_1,\cdots, y_m\}\cup U)$ must be an ideal of $S$. If
further $a^2\not=0$, then $S\setminus \{y_1,\cdots, y_m\}$ is a
sub-semigroup of $S$.  Also by \cite[Theorem 2.1]{WuLu}, $S\setminus
(\{y_1,\cdots, y_m\}\cup U\cup V)$ is a sub-semigroup of $S\setminus
(\{y_1,\cdots, y_m\}\cup U)$, and thus a sub-semigroup of $S$.

For $m=2, n=2$, $U=\{u\}$ and $V=\{v,\ol{v}\}$, it is not very hard
to construct a semigroup $T$ such that $\G(T)=G-\{y_1,y_2\}$
following the way mentioned above. Then after a rather complicated
calculation, we succeed in adding two vertices $y_1,y_2$ to this
table such that $\G(S)=G$. The multiplication on $S$ is listed in
Table 3 and the detailed verification for the associativity is
omitted here:

$${\text{ Table 3}}$$
$$\begin{array}{c|cccccccccc}
\cdot& a&d&b&x_1&x_2&y_1&y_2&u&v&\overline{v}\\
\hline a&a&0&0&0&0&0&0&0&a&a\\
d&0&d&0&0&0&d&d&d&0&0\\
b&0&0&0&b&b&0&0&b&b&b\\
x_1&0&0&b&x_1&x_2&b&b&x_1&x_1&x_1\\
x_2&0&0&b&x_2&x_2&b&b&x_2&x_2&x_2\\
y_1&0&d&0&b&b&d&d&y_1&b&b\\
y_2&0&d&0&b&b&d&d&y_2&b&b\\
u&0&d&b&x_1&x_2&y_1&y_2&u&x_1&x_1\\
v&a&0&b&x_1&x_2&b&b&x_1&v&v\\
\overline{v}&a&0&b&x_1&x_2&b&b&x_1&v&v\\
\end{array}
$$

Notice that $S\setminus\{x_1,x_2\}$ is not a sub-semigroup of $S$
since $uv=x_1$. Notice also that $S\setminus (U\cup
\{x_1,x_2,\cdots, x_n\})$ is a sub-semigroup of $S$.

 We remark that the construction in
Table 3 can be routinely extended for all $n\ge 1,m\ge 1, |U|\ge 0$
and $|V|\ge 0$, where each of $m,n, |U|,|V|$ could be a finite or an
infinite cardinal number. In other words, each graph in Fig.3 has a
corresponding semigroup for any finite or infinite $n\ge 1,m\ge 1,
|U|\ge 0$ and $|V|\ge 0$.

\vs{3mm}\nin{\bf Remark 3.2.} {\it Consider the graph $G$ in Fig.3
and assume that $n\ge 1,$ $m\ge 1, |U|\ge 0,$$|V|\ge 1$.}

(1) {\it If we add an  end vertex $w$ which is adjacent to $b$, then
the resulting graph $\ol{G}$ has no corresponding zero-divisor
semigroup, even if  $U=\eset$.}

(2) {\it If we add a vertex $w$ such that $N(w)=\{b,d\}$, then the
resulting graph $H$ has no corresponding zero-divisor semigroup,
even if  $U=\emptyset$.}

\vs{2mm}\nin{\bf Proof.} (1) Assume $v\in V$. We only need consider
the case when $U=\emptyset$. Suppose that $\ol{G}$ is the
zero-divisor graph of a semigroup $S$ with $V[\G(S)]=V(\ol{G})$. By
Proposition 2.2(2), we have $bx_1=bv=b$ and $dy_1=dw=d$. Clearly,
$aw, av\in Ann(x_1)\cap Ann(y_1)=\{a,0\}$, and thus $aw=a$ and
$av=a$. As $awv=av=a$, we have $wv\in [Ann(b)\cap Ann(d)]\setminus
Ann(a)$. That means $wv=a$ and $a^2\neq 0$. We have $y_1wv=y_1a=0$
and $x_1wv=x_1a=0$. Thus $y_1w=x_1w=d$ and $y_1v=x_1v=b$ by Lemma
2.1. Consider $x_1y_1v$. We have
$b=x_1b=x_1(y_1v)=y_1(x_1v)=y_1b=0$, a contradiction. The
contradiction shows that $\ol{G}$ has no corresponding semigroup.

(2) Assume $v\in V$. We only need consider the case $U=\emptyset$.
Suppose that $H$ is the zero-divisor graph of a semigroup $S$. We
have $bS\subseteq Ann(y)\cap Ann(w)\subseteq \{0,b\}$. Clearly,
$wv\in [Ann(d)\cap Ann(b)]\setminus Ann(a)\subseteq \{a,w\}$ since
$wva=wa=a$. If $wv=a$, then we have $wvx_1=ax_1=0$ and
$wvy_1=ay_1=0$, which means $wx_1, wy_1\in Ann(v)\subseteq
\{0,v,d\}$. As $wx_1,wy_1\in Ann(a)\setminus \{0\}$, we have
$wx_1=wy_1=d$. Then $0=dx_1=wy_1x_1=y_1d=d$, a contradiction. Now
assume $wv=w$ and consider $wx_1$. $wx_1\in Ann(a)\cap Ann(d)\cap
Ann(b)\subseteq \{0,a,b,d\}$. We claim $wx_1\neq d$  since
otherwise, $d=wx_1=wv\cdot x_1=wx_1\cdot v=dv=0$, a contradiction.
In a similar way we prove $wy_1\in \{a,b\}$. Moreover, $wy_1\cdot
x_1=y_1\cdot wx_1=0$ whether $wx_1=a$ or $wx_1=b$. Thus $wy_1=a$. As
$x_1^2w=y_1^2w=x_1y_1w=0$, we have $x_1y_1, x_1^2, y_1^2 \in
Ann(a)\cap Ann(w)\subseteq \{b,d,0\}$, but $x_1y_1\neq 0$.
 Now consider $x_1y_1$. We conclude $x_1y_1=d$ since otherwise, $x_1y_1=b$
and it implies $b=bx_1=x_1y_1x_1=x_1^2y_1\neq b$, a contradiction.
Finally, $x_1y_1=d$ implies $d=dy_1=y_1x_1y_1=x_1y_1^2\in \{0,b\}$,
a contradiction. This completes the proof. \quad $\Box$

\vs{3mm}Now come back to the structure of semigroup graphs $G$
satisfying the main assumption in Theorem 2.4. We use notations used
in its proof. The vertex set of the graph was decomposed into four
mutually disjoint nonempty parts, i.e., $V(G)= \{a,b\}\cup
C(a,b)\cup B\cup L$, where after taking a $c$ in $C(a,b)$
$$B=\{v\in V(G)\,|\, v\notin C(a,b),\,d(c,v)=2\},\,L=\{v\in V(G)\,|\, d(c,v)=3\}.$$
(For example, for the graph $G$ in Fig.3,  $C(a,b)=\{y_j\}$,
$B=U\cup \{d\}\cup \{x_i\},\,L=V$. In particular, $L$ consists of
end vertices.) By \cite[Theorem 1(4)]{DD}, for each pair $x,y$ of
nonadjacent vertices of $G$, there is a vertex $z$ with $N(x)\cup
 N(y)\seq \ol{ N(z)}$. Then we have the following observations:

 (1) {\it No two vertices in $L$ are adjacent in $G$.} Thus a vertex of
 $L$ is either an end vertex or is adjacent to at least two vertices
 in $B$. In particular, the subgraph induced on $L$ is a completely
 discrete graph.

 (2) A vertex in $B$ is adjacent to either $a$ or $b$. {\it If a vertex $k$ in $B$ is adjacent to a vertex $l$ in $L$,
 then $k$ is adjacent to both $a$ and $b$.} Thus $B$ consists of four
 parts: end vertices in $T_a$ that are adjacent to $a$, end vertices in $T_b$ that are adjacent to
 $b$, vertices in $B_2$ that are adjacent to both $a$ and $b$, and vertices in $B_1$ that are adjacent
 to one of $a,b$ and at the same time adjacent to another vertex in
 $B$. By Example 3.1, the structure of the induced subgraph on $B_1\cup B_2$ seems to be
 complicated. In the following, we will give a complete
 classification of the semigroup graphs $G$ with $|B_1\cup B_2|\le
 2$.

 \vs{2mm} First, consider the case $|B_1\cup B_2|=1$.

 \vs{3mm}\nin{\bf Theorem 3.3.} {\it Let $G$ be a connected,
simple graph with diameter 3. Assume that there exist two adjacent
vertices $a,b$ in $V(G)$, and assume that there exist a vertex
$c_1\in C(a,b)$ and a vertex $w\in V(G)$ such that $d(c_1,w)=3$.
Assume further that $|B\setminus(T_a\cup T_b)|=1$. Then $G$ is a
semigroup graph if and only if the following conditions hold:} (1)
{\it $1\le |C(a,b)|\leq \infty, 1\le |W|\leq \infty$ and $W$
consists of end vertices, where $W=\{s\in V(G)\,|\, d(c_1, s)=3\}$.}
(2) {\it either $T_a=\eset$ or $T_b=\eset.$ (see Fig.4 with
$V=\eset$.)}

\vs{2mm}\nin{\bf Proof.} As $|B\setminus(T_a\cup T_b)|=1$,
$B\setminus(T_a\cup T_b)=B_2$. By the previous observations, we need
only prove the following two facts.

(1)  If $|T_a|\geq 0$ and $T_b=\emptyset$, then $G$ is a subgraph of
Fig.3 with $C(a,d)=\emptyset$. (see also Fig.4 with $V=\eset$.) We
claim that $G$ is a semigroup graph. In fact, if $U=\eset$, delete
the three rows and the three columns involving $x_1,x_2$ and $u$ in
Table 3 to obtain an associative multiplication on $S_1=S\setminus
(U\cup \{x_1,x_2,\cdots, x_n\})$. Clearly, $\G(S_1)=G$ for
$|C(a,b)|=2=|V|, |C(a,d)|=0=|U|$ in Fig.3. Also, the table can be
extended for any finite or infinite $|C(a,b)|\ge 1$ and $|V|\ge 0$
while $|U|= 0$. If $|U|>0$, then we work out a corresponding
associative multiplication table listed in Table 4, for
$C(a,b)=\{y_1,y_2\}$, $U=\{u_1,u_2\}$, $V=\{v_1,v_2\}$.

$${\text{ Table 4}}$$
$$\begin{array}{c|ccccccccc}
\cdot& a&b&d&y_1&y_2&u_1&u_2&v_1&v_2\\
\hline
a&0&0&0&0&0&0&0&a&a\\
b&0&a&0&0&0&a&a&b&b\\
d&0&0&d&d&d&d&d&0&0\\
y_1&0&0&d&d&d&d&d&a&a\\
y_2&0&0&d&d&d&d&d&a&a\\
u_1&0&a&d&d&d&y_1&y_1&b&b\\
u_2&0&a&d&d&d&y_1&y_1&b&b\\
v_1&a&b&0&a&a&b&b&v_1&v_1\\
v_2&a&b&0&a&a&b&b&v_1&v_1\\
\end{array}$$

Clearly, the table can be extended for all finite or infinite
$|C(a,b)|\ge 1$, $|U|\ge 1$ and $|V|\ge 1$. This completes the
proof.

\vs{2mm} (2) If both $|T_a|>0$ and $|T_b|>0$, then we conclude that
$G$ is not a semigroup graph.

\begin{center}
\setlength{\unitlength}{0.13cm}
\begin{picture}(74,42)
\thinlines
\drawdot{36.0}{6.0}
\drawdot{52.0}{14.0}
\drawdot{20.0}{14.0}
\drawdot{10.0}{14.0}
\drawdot{62.0}{14.0}
\drawdot{48.0}{6.0}
\drawpath{36.0}{6.0}{48.0}{6.0}
\drawpath{36.0}{6.0}{52.0}{14.0}
\drawpath{52.0}{14.0}{62.0}{14.0}
\drawpath{36.0}{6.0}{20.0}{14.0}
\drawpath{20.0}{14.0}{52.0}{14.0}
\drawpath{10.0}{14.0}{20.0}{14.0}
\drawdot{24.0}{24.0}
\drawdot{28.0}{24.0}
\drawdot{32.0}{24.0}
\drawdot{36.0}{24.0}
\drawdot{40.0}{24.0}
\drawdot{48.0}{24.0}
\drawpath{20.0}{14.0}{24.0}{24.0}
\drawpath{24.0}{24.0}{52.0}{14.0}
\drawpath{52.0}{14.0}{28.0}{24.0}
\drawpath{28.0}{24.0}{20.0}{14.0}
\drawpath{20.0}{14.0}{48.0}{24.0}
\drawpath{48.0}{24.0}{52.0}{14.0}
\drawcenteredtext{34.0}{36.0}{\tiny Fig.4.}
\drawcenteredtext{6.0}{14.0}{$U$}
\drawcenteredtext{66.0}{14.0}{$V$}
\drawcenteredtext{52.0}{6.0}{$W$}
\drawcenteredtext{30.0}{6.0}{$d$}
\drawcenteredtext{18.0}{12.0}{$a$}
\drawcenteredtext{54.0}{12.0}{$b$}
\drawcenteredtext{20.0}{28.0}{$c_1$}
\drawcenteredtext{28.0}{28.0}{$c_2$}
\drawcenteredtext{48.0}{28.0}{$c_n$}
\drawpath{56.0}{6.0}{56.0}{6.0}
\end{picture}\par
\end{center}

In fact, in this case, $G$ is a graph in Fig.4, where $|W|\ge 1,
|U|\ge 1,|V|\ge 1$. Assume $u\in U$, $v\in V$, $w\in W$ and $c\in
C(a,b)$. We now proceed to prove that such a graph does not have a
corresponding semigroup.

Suppose that $G$ is the zero-divisor graph of a semigroup $S$ with
$V[\G(S)]=V(G)$. By Proposition 2.2(2), we have $du=dv=dc_1=d$. Then
$uc_1d=ud=d$, which implies $uc_1\in [Ann(a)\cap Ann(b)]\setminus
Ann(d)\subseteq \{c_i,d\}$.

Assume $uc_1=d$. Then $uc_1w=dw=0$, and thus $c_1w=a$ by Lemma 2.1.
As $c_1wv=av=a$, we have $wv\in [Ann(d)\cap Ann(b)]\setminus
Ann(c_1)\subseteq \{d\}$, thus $wv=d$. Then $a=wvc_1=dc_1=d$, a
contradiction.

So $uc_1=c_i$, and therefore $wuc_1=wc_i\neq 0$. We have $$wu\in
[Ann(a)\cap Ann(d)]\setminus Ann(c_1)\subseteq \{d\},$$and thus
$wu=d$. Then $b=bw=buw=bd=0$, a contradiction. This completes the
proof. \quad $\Box$

\vs{2mm}A natural question arising from Example 3.1 is if $L$ only
consists of end vertices. The following example shows this is not
the case.

 \vs{2mm}\nin{\bf Example 3.4.} {\it  Consider the graph $G$ in Fig.5,
 where $C(a,b)=\{c_1,c_2,\cdots,c_m \}$, $L=\{y_1,y_2,\cdots, y_n\}$, $B=\{x_1, x_2 \}\cup V$ ($m\geq 1, n\ge 1, |V|\ge 0$ ) and $V$ consists of end vertices adjacent to
 $b$. Notice that each of $m,n$ and $|V|$ could be finite or
 infinite. We conclude that each graph in Fig.5 has a corresponding
 zero-divisor semigroup.}

\begin{center} 
\setlength{\unitlength}{0.13cm}
\begin{picture}(68,44)
\thinlines
\drawdot{58.0}{20.0}
\drawpath{12.0}{20.0}{16.0}{30.0}
\drawpath{16.0}{30.0}{46.0}{20.0}
\drawpath{46.0}{20.0}{22.0}{30.0}
\drawpath{22.0}{30.0}{12.0}{20.0}
\drawpath{12.0}{20.0}{40.0}{30.0}
\drawpath{40.0}{30.0}{46.0}{20.0}
\drawpath{12.0}{20.0}{12.0}{14.0}
\drawpath{12.0}{14.0}{46.0}{20.0}
\drawpath{46.0}{20.0}{46.0}{14.0}
\drawpath{46.0}{14.0}{12.0}{20.0}
\drawpath{12.0}{14.0}{16.0}{8.0}
\drawpath{16.0}{8.0}{46.0}{14.0}
\drawpath{46.0}{14.0}{22.0}{8.0}
\drawpath{22.0}{8.0}{12.0}{14.0}
\drawpath{12.0}{14.0}{40.0}{8.0}
\drawpath{40.0}{8.0}{46.0}{14.0}
\drawdot{12.0}{20.0}
\drawdot{12.0}{14.0}
\drawpath{18.0}{34.0}{18.0}{34.0}
\drawcenteredtext{16.0}{32.0}{$c_1$}
\drawcenteredtext{40.0}{32.0}{$c_m$}
\drawdot{46.0}{20.0}
\drawdot{46.0}{14.0}
\drawpath{12.0}{20.0}{46.0}{20.0}
\drawcenteredtext{8.0}{20.0}{$a$}
\drawcenteredtext{48.0}{22.0}{$b$}
\drawcenteredtext{22.0}{32.0}{$c_2$}
\drawcenteredtext{32.0}{38.0}{\tiny Fig.5.}
\drawcenteredtext{60.0}{20.0}{$V$}
\drawcenteredtext{8.0}{14.0}{$x_1$}
\drawcenteredtext{50.0}{14.0}{$x_2$}
\drawcenteredtext{16.0}{6.0}{$y_1$}
\drawcenteredtext{22.0}{6.0}{$y_2$}
\drawcenteredtext{42.0}{6.0}{$y_n$}
\drawdot{40.0}{30.0}
\drawdot{16.0}{30.0}
\drawdot{22.0}{30.0}
\drawdot{26.0}{30.0}
\drawdot{30.0}{30.0}
\drawdot{34.0}{30.0}
\drawpath{46.0}{20.0}{58.0}{20.0}
\drawdot{16.0}{8.0}
\drawdot{22.0}{8.0}
\drawdot{26.0}{8.0}
\drawdot{30.0}{8.0}
\drawdot{34.0}{8.0}
\drawdot{40.0}{8.0}
\end{picture}\par
\end{center}

\Proof We need only  work out a corresponding associative
multiplication table for $|V|=m=n=2$. We use Theorem 2.4 and list
the associative multiplication in Table 5. Clearly, the table can be
extended for all finite or infinite $m,n\ge 1$, and $|V|\ge 0$.
$${\text{ Table 5}}$$
$$\begin{array}{c|cccccccccc}
\cdot& a&b&c_1&c_2&v_1&v_2&x_1&x_2&y_1&y_2\\
\hline
a&a&0&0&0&a&a&0&0&a&a \\
b&0&0&0&0&0&0&0&0&b&b\\
c_1&0&0&x_1&x_1&x_1&x_1&x_1&x_1&b&b\\
c_2&0&0&x_1&x_1&x_1&x_1&x_1&x_1&b&b\\
v_1&a&0&x_1&x_1&v_1&v_1&x_1&x_1&a&a\\
v_2&a&0&x_1&x_1&v_1&v_1&x_1&x_1&a&a\\
x_1&0&0&x_1&x_1&x_1&x_1&x_1&x_1&0&0\\
x_2&0&0&x_1&x_1&x_1&x_1&x_1&x_2&0&0\\
y_1&a&b&b&b&a&a&0&0&y_1&y_1\\
y_2&a&b&b&b&a&a&0&0&y_1&y_1\\
\end{array}$$ This completes the proof.\quad $\Box$

\vs{2mm} We have three remarks to Example 3.4.

\vs{2mm}(1) {\it Let $n\ge 1, m\ge 1$. If we add to $G$ in Fig.5 an
end vertex $u$ such that $au=0$, then the resulting graph $\ol{G}$
has no corresponding zero-divisor semigroup.}

\vs{2mm}\nin{\bf Proof.} (1) Suppose to the contrary  that $\ol{G}$
is the zero-divisor graph of a semigroup $P$ with
$V[\G(P)]=V(\ol{G})$. By Proposition 2.2(2), we have $a^2\in
\{0,a\}$ and $b^2\in \{0,b\}$. First, we have $v_1y_1 \in Ann(b)\cap
Ann(x_1)\cap Ann(x_2)=\{a,b,0\}$ and similarly, $uy_1, c_1y_1\in
\{a,b\}$. Then $v_1y_1=a$ and $a^2=a$ since $av_1y_1=ay_1=a$. On the
other hand, $auy_1=0$ and it implies $uy_1=b$. Similarly, we have
$c_1y_1=b$. Consider $c_1uy_1$. We have
$b=ub=u(c_1y_1)=c_1(uy_1)=c_1b=0$, a contradiction. This completes
the proof. \quad $\Box$

\vs{2mm}(2) {\it Let $n\ge 1, m\ge 1$. If we add to $G$ in Fig.5  an
end vertex $y$ such that $yx_1=0$, then the resulting graph $\ol{G}$
has no corresponding zero-divisor semigroup, whether or not $T_b=
\emptyset $.}

\vs{2mm}\nin{\bf Proof.} (2) Assume $\{y_1,y\}\seq L$, where $y$ is
an end vertex adjacent to $x_1$. Suppose to the contrary  that
$\ol{G}$ is the zero-divisor graph of a semigroup $P$ with
$V[\G(P)]=V(\ol{G})$. First, $x_2y\in Ann(a)\cap Ann(x_1)\cap
Ann(y_1)=\{x_1,0\}$. Thus $x_2y=x_1$, and hence $x_1^2=0$. By
Proposition 2.2(2), we have $c_1x_1=x_1$ and therefore,
$c_1^2x_1=x_1$. Thus $c_1^2\in \{c_i,x_2\,|\,i\}$. We have
$c_1^2y_1=0$ since $c_1y_1\in Ann(a)\cap Ann(x_1)\cap
Ann(x_2)=\{a,b,0\}$. Since $c_1^2=x_2$, $c_1y\in \{a,b,x_1\}$ and
$c_1^2y=x_2y=x_1$, it follows that $c_1y=x_1$. Finally,
$c_1yy_1=x_1y_1=0$ and by Lemma 2.1, we have $c_1y_1=x_1$,
contradicting $c_1y_1\in \{a,b\}$. This completes the proof.\quad
$\Box$

\vs{2mm}(3) {\it Let $n\ge 1, m\ge 1$ and assume $V=\eset$ in Fig.5.
If further we add to $G$  an edge connecting $x_1$ and $x_2$, then
the resulting graph $\ol{G}$ has no corresponding zero-divisor
semigroup.}

\vs{2mm}\nin{\bf Proof.} Suppose to the contrary that $\ol{G}$ is
the zero-divisor graph of a semigroup $P$ with $V[\G(P)]=V(\ol{G})$.
By Lemma 2.1, we have $c_1x_1\in Ann(y_1)=\{x_1,x_2,0\}$ and
similarly, $c_1x_2\in \{x_1,x_2\}$. Then we have $c_1^2x_1\neq 0$
and $c_1^2x_2\neq 0$, which means $c_1^2 \in [Ann(a)\cap
Ann(b)]\setminus [Ann(x_1)\cup Ann(x_2)]=\{c_i|i=1,2,\cdots,m\}$
since $c_1^2\neq 0$ by Lemma 2.1. Similarly, we have $y_1^2\in
\{y_i|i=1,2,\cdots,n\}$. Clearly, we have $c_1y_1\in Ann(a)\cap
Ann(x_1)\subseteq \{a,b,x_1,x_2,0\}$. Then as $c_1^2y_1=c_iy_1\neq
0$ for some $i\in \{1,2,\cdots, m\}$, we have $c_1y_1\in
\{x_1,x_2\}$. Finally, $0=(c_1y_1)y_1=c_1y_1^2=c_1y_i\neq 0$ (for
some $i\in \{1,2,\cdots, n\}$), a contradiction. This completes the
proof.\quad $\Box$

\vs{3mm}Combining the above results, we now classify all semigroup
graphs satisfying the main assumption of Theorem 2.4 with $|B_1\cup
B_2|=2$:

\vs{3mm}\nin{\bf Theorem 3.5.} {\it Let $G$ be a connected, simple
graph with diameter 3. Assume that there exist two adjacent vertices
$a,b$ in $V(G)$, and assume that there exist a vertex $y_1\in V(G)$
and a vertex $c_1\in C(a,b)$ such that $d(c_1,y_1)=3$. Assume
further $|B\setminus(T_a\cup T_b)|=2$. }

(1) {\it If $B_2=B\setminus(T_a\cup T_b)$, then $G$ is  a semigroup
graph if and only if $G$ is a  graph in Fig.5, where $1\le m\le
\infty,1\le n\le \infty$ and $0\le |V|\le \infty.$}

(2) {\it If $|B_2|=1$,  then $G$ is  a semigroup graph if and only
$G$ is a  graph in Fig.3, where $n=1,1\le m\le \infty,$ $1\le |V|\le
\infty$, $0\le |U|\le \infty.$}

\Proof (1) By Example 3.4, each graph in Fig.5 is a semigroup graph.
Clearly, $B_2=B\setminus(T_a\cup T_b)$ and it consists of two
vertices. Conversely, the result follows from \cite[Theorem
2.1]{WuLu} and the three remarks after Example 3.4.

(2) If $|B_2|=1$, then assume $B\setminus(T_a\cup T_b)=\{x_1,x_2\}$,
where $a-x_2-b$. In this case, $x_1-x_2$ in $G$. If $x_1-a$ in $G$,
then there is no end vertex adjacent to $x_1$. In this subcase, $G$
is a semigroup graph if and only if $T_b=\eset$ by Example 3.1 and
Remark 3.2(1), the case of $|C(a,d)|=1$. The other subcase is
$x_1-b$ in $G$, and it is the same with the above subcase. This
completes the proof. \quad $\Box$

\vs{3mm}It is natural to ask the following {\it question}: Can one
give a complete classification of  semigroup graphs $G=\G(S)$ with
$|B_1\cup B_2=n|$ for any $n\ge 3$? At present, it seems to be a
rather difficult question.

\vs{3mm}Add two end vertices to two vertices of the complete graph
$K_n$ to obtain a new graph, and denote the new graph as $K_n +2$.
By \cite[Theorem 2.1]{WuLu2}, $K_n +2$ has a unique zero-divisor
semigroup $S$ such that $\G(S)\cong K_n+2$ for each $n\ge 4$. Having
Theorem 2.4 in mind, it is natural to consider graphs obtained by
adding some caps to $K_n+2$.

\vs{2mm}\nin{\bf Example 3.6.} {\it Consider the graph $G$ in Fig.6.
The subgraph $G_1$ induced on the vertex subset
$S^*=\{a,b,x_1,x_2,y_1,y_2\}$ is the graph $K_4+2$, i.e., $K_4$
together with two end vertices $y_1,y_2$. Then $G_1$ has a unique
corresponding zero-divisor semigroup $S=S^*\cup \{0\}$ by
\cite[Theorem 2.1]{WuLu2}. We can work out the corresponding
associative multiplication table, and list it in Table 6:}

\begin{center}
 \setlength{\unitlength}{0.13cm}
\begin{picture}(70,58)
\thinlines
\drawpath{24.0}{42.0}{44.0}{42.0}
\drawpath{44.0}{42.0}{44.0}{24.0}
\drawpath{44.0}{24.0}{24.0}{24.0}
\drawpath{24.0}{24.0}{24.0}{42.0}
\drawpath{24.0}{24.0}{8.0}{24.0}
\drawpath{44.0}{24.0}{60.0}{24.0}
\drawdot{24.0}{24.0}
\drawdot{24.0}{42.0}
\drawdot{44.0}{42.0}
\drawdot{44.0}{24.0}
\drawdot{60.0}{24.0}
\drawdot{8.0}{24.0}
\drawpath{24.0}{42.0}{44.0}{24.0}
\drawpath{24.0}{24.0}{44.0}{42.0}
\drawcenteredtext{20.0}{42.0}{$a$}
\drawcenteredtext{48.0}{42.0}{$b$}
\drawcenteredtext{20.0}{22.0}{$x_1$}
\drawcenteredtext{48.0}{22.0}{$x_2$}
\drawcenteredtext{8.0}{22.0}{$y_1$}
\drawcenteredtext{60.0}{22.0}{$y_2$}
\drawcenteredtext{34.0}{52.0}{\tiny Fig.6.}
\drawdot{34.0}{20.0}
\drawdot{34.0}{16.0}
\drawdot{34.0}{8.0}
\drawdot{34.0}{12.0}
\drawpath{24.0}{24.0}{34.0}{20.0}
\drawpath{34.0}{20.0}{44.0}{24.0}
\drawpath{24.0}{24.0}{34.0}{16.0}
\drawpath{34.0}{16.0}{44.0}{24.0}
\drawpath{24.0}{24.0}{34.0}{8.0}
\drawpath{34.0}{8.0}{44.0}{24.0}
\drawcenteredtext{34.0}{22.0}{$c_1$}
\drawcenteredtext{34.0}{18.0}{$c_2$}
\drawcenteredtext{34.0}{6.0}{$c_n$}
\drawdot{34.0}{46.0}
\drawpath{24.0}{42.0}{34.0}{46.0}
\drawpath{34.0}{46.0}{44.0}{42.0}
\drawcenteredtext{34.0}{48.0}{$c$}
\drawcenteredtext{16.0}{32.0}{$d$}
\drawdot{20.0}{32.0}
\drawpath{24.0}{42.0}{20.0}{32.0}
\drawpath{20.0}{32.0}{24.0}{24.0}
\drawpath{24.0}{24.0}{24.0}{24.0}
\drawpath{52.0}{32.0}{52.0}{32.0}
\end{picture}\par
\end{center}

\nin

$${\text{ Table 6}}$$
$$\begin{array}{c|cccccc}
\cdot& a&b&x_1&x_2&y_1&y_2\\
\hline
a&a&0&0&0&a&a \\
b&0&0&0&0&x_2&x_1\\
x_1&0&0&0&0&0&x_1\\
x_2&0&0&0&0&x_2&0\\
y_1&a&x_2&0&x_2&y_1&a\\
y_2&a&x_1&x_1&0&a&y_2\\
\end{array}$$

(1) {\it If we add to $G_1$ a vertex $c$ such that $N(c)=\{a,b\}$,
then the resulting graph $H_1$ has no corresponding zero-divisor
semigroup.}

(2) {\it If we add to $G_1$ a vertex $d$ such that $N(d)=\{a,x_1\}$,
then the resulting graph $H_2$ has no corresponding zero-divisor
semigroup.}

(3) {\it If we add to $G_1$ vertices $c_i$ ($i\in I$) such that
$N(c_i)=\{x_1,x_2\}$, then the resulting graph $H$ has corresponding
zero-divisor semigroups, where $I$ could be any finite or infinite
index set.}

In each of the above three cases, we say that {\it  a cap is added
to the subgraph $K_4+2$}.

\vs{2mm}\nin{\bf Proof.} (1) Suppose that $H_1$ is the zero-divisor
graph of a semigroup $S_1$ with $V[\G(S_1)]=V(H_1)$. Then by Theorem
2.4, $S$  is an ideal of $S_1=S\cup \{c\}$. Thus we only need check
the associative multiplication of $S_1$ based on the table of $S$
already given in Table 6. First, we have $cx_2=x_2$ by Proposition
2.2(2). Consider $y_1bc$. Clearly,
$0=0y_1=(cb)y_1=c(by_1)=cx_2=x_2$, a contradiction. This completes
the proof.

\vs{2mm} (2) Suppose that $H_2$ is the zero-divisor graph of a
semigroup $S_2=S\cup \{d\}$ with $V[\G(S_2)]=V(H_2)$. If
$x_1^2\not=0$, then by Theorem 2.4(2), $S$ is a sub-semigroup of
$S_2$. Then $\G(S)=K_4+2$, and it implies $x_1^2=0$ by Table 6, a
contradiction. In the following we assume $x_1^2=0$.

By Lemma 2.1, we have $d^2\neq 0$, and thus $ay_1,ay_2\in
Ann(d)=\{a,x_1,0\}$. Clearly $ay_1\neq 0$ and we can have $ay_2=a$.
(Otherwise, $ay_2=x_1$ and we have $0=x_1y_1=ay_2y_1=(ay_1)y_2\neq
0$, a contradiction.) Then $ay_1y_2\neq 0$, and thus $y_1y_2\in
[Ann(x_1)\cap Ann(x_2)]\setminus Ann(a)$. It means $y_1y_2=a$ and
$a^2\neq 0$. Clearly $by_1y_2=0$, and thus $by_1=x_2$, $by_2=x_1$ by
Lemma 2.1. Similarly, $cy_1y_2=ca=0$  and thus $cy_1=x_2$,
$cy_2=x_1$. Finally, consider $bcy_1$. We have
$0=bx_2=b(cy_1)=c(by_1)=cx_2=x_2$, a contradiction. This completes
the proof.

\vs{2mm} (3) Suppose that $H$ is the subgraph of $G$ in Fig.6
induced on the vertex set $S^*\cup \{c_i\,|\, i\in I\}$. Assume that
$H$ is the zero-divisor graph of a semigroup $P$ with
$V[\G(P)]=V(H)$. Clearly, it dose not satisfy the condition of
Theorem 2.4. For $|I|=2$, we work out an associative multiplication
table and list it in Table 7:

$${\text{ Table 7}}$$
$$\begin{array}{c|cccccccc}
\cdot& a&b&x_1&x_2&y_1&y_2&c_1&c_2\\
\hline
a&a&0&0&0&a&a&a&a \\
b&0&b&0&0&b&b&b&b\\
x_1&0&0&x_1&0&0&x_1&0&0\\
x_2&0&0&0&x_2&x_2&0&0&0\\
y_1&a&b&0&x_2&y_1&c_1&c_1&c_1\\
y_2&a&b&x_1&0&c_1&y_2&c_1&c_1\\
c_1&a&b&0&0&c_1&c_1&c_1&c_1\\
c_2&a&b&0&0&c_1&c_1&c_1&c_1
\end{array}$$
\nin The table can be easily extended for any finite or infinite
index set $I$. \quad $\Box$

\vs{3mm} We remark that in Example 3.6, replace $K_4$ by $K_n$ for
any $n\ge 5$,  the results still hold. There exists no difficulty to
generalize the proofs to the general cases. Thus we have proved the
following general result.

\vs{3mm}\nin{\bf Theorem 3.7.} {\it Assume $n\ge 4$ and let
$G=K_n+2$ be the complete graph $K_n$ together with two end
vertices. Add some (finite or infinite) caps to the subgraph $K_n$
to obtain a new graph $H$ such that $G$ is a subgraph of $H$. Then
$H$ is a semigroup graph if and only if each of the gluing vertices
is adjacent to an end vertex in $G$. }

\vs{4mm}

\end{document}